\documentclass[titlepage,draft,12pt]{article} 
\usepackage{amsfonts,latexsym}
\usepackage[a4paper]{geometry}

\date{}

\date{\mbox{}\bigskip\\
\emph{Dedicated to Professor Sergio Spagnolo on the occasion
 of his $70$th birthday}%\bigskip\\
}

\newcommand{\ep}{\varepsilon}
\newcommand{\qed}{{\penalty 10000\mbox{$\quad\Box$}}}
\newcommand{\re}{\mathbb{R}}
\newcommand{\C}{\mathbb{C}}
\newcommand{\n}{\mathbb{N}}
\newcommand{\aholder}{$\alpha$-H\"{o}lder}
\newcommand{\hold}{\mbox{\textup{H\"{o}ld}}_{\alpha}}

\newtheorem{thm}{Theorem}[section]
\newtheorem{thmbibl}{Theorem}

\newtheorem{prop}[thm]{Proposition}

\newtheorem{ex}[thm]{Example}
\newtheorem{lemma}[thm]{Lemma}

\title{Higher order Glaeser inequalities and optimal regularity of
roots of real functions}

\author{Marina Ghisi\vspace{1ex}\\ 
{\normalsize Universit\`a degli Studi di Pisa} \\
{\normalsize Dipartimento di Matematica ``Leonida Tonelli''}\\ 
{\normalsize PISA (Italy)}\\
{\normalsize e-mail: \texttt{ghisi@dm.unipi.it}}
\and
Massimo Gobbino\vspace{1ex}\\ 
{\normalsize Universit\`a degli Studi di Pisa} \\
{\normalsize Dipartimento di Matematica Applicata ``Ulisse Dini''}\\ 
{\normalsize PISA (Italy)}\\  
{\normalsize e-mail: \texttt{m.gobbino@dma.unipi.it}}}

\begin{document}
\maketitle
\begin{abstract}
	We prove a higher order generalization of Glaeser inequality,
	according to which one can estimate the first derivative of a
	function in terms of the function itself, and the H\"{o}lder
	constant of its $k$-th derivative.
	
	We apply these inequalities in order to obtain \emph{pointwise}
	estimates on the derivative of the $(k+\alpha)$-th root of a
	function of class $C^{k}$ whose derivative of order $k$ is
	\aholder\ continuous.  Thanks to such estimates, we prove that the
	root is not just absolutely continuous, but its derivative has a
	higher summability exponent.
	
	Some examples show that our results are optimal.
	
\vspace{1cm}

\noindent{\bf Mathematics Subject Classification 2000 (MSC2000):}
26A46, 26B30, 26A27.

\vspace{1cm} 

\end{abstract}

%%%%%%%%%%%%%%%%%%%%%
%                   %
%   Inizio lavoro   %
%                   %
%%%%%%%%%%%%%%%%%%%%%
 
\section{Introduction}

Let $v\in C^{2}(\re)$ be a function such that either $v(x)\geq 0$ or
$v(x)\leq 0$ for every $x\in\re$.  Let us assume that $v''(x)$ is
bounded.  Then one has that
\begin{equation}
	|v'(x)|^{2}\leq 2|v(x)|\cdot\sup_{x\in\re}|v''(x)|
	\quad\quad
	\forall x\in\re.
	\label{glaeser}
\end{equation}

This is known as Glaeser inequality (see~\cite{glaeser}).  If we set
$u(x):=|v(x)|^{1/2}$, then estimate (\ref{glaeser}) implies that
$$|u'(x)|=\frac{|v'(x)|}{2|v(x)|^{1/2}}\leq
\left\{\frac{1}{2}\sup_{x\in\re}|v''(x)|\right\}^{1/2}$$
for every $x\in\re$ such that $v(x)\neq 0$.  In other words, a bound
on the second derivative of some function $v(x)$ with constant sign
yields a bound on the first derivative of the square root of $|v(x)|$.

The first question we address in this paper is how to obtain bounds 
on the first derivative of higher order roots of $v$. We prove that 
for every integer $k\geq 2$ there exists a constant $C(k)$, depending 
only upon $k$, such that
\begin{equation}
	|v'(x)|^{k+1}\leq C(k)\cdot|v(x)|^{k}\cdot\sup_{x\in\re}|v^{(k+1)}(x)|
	\quad\quad
	\forall x\in\re
	\label{glaeser-k}
\end{equation}
for every $v\in C^{k+1}(\re)$ such that either $v(x)\geq 0$ or
$v(x)\leq 0$ for every $x\in\re$, and either $v'(x)\geq 0$ or
$v'(x)\leq 0$ for every $x\in\re$.  This means that we can estimate
the first derivative of $u(x):=|v(x)|^{1/(k+1)}$ (even without
absolute value, if $k$ is even) in terms of the $(k+1)$-th derivative
of $v$, namely 
$$|u'(x)|=\frac{1}{k+1}\cdot\frac{|v'(x)|}{|v(x)|^{k/(k+1)}}\leq
\frac{1}{k+1}\left\{C(k)\sup_{x\in\re}|v^{(k+1)}(x)|\right\}^{1/(k+1)}$$
for every $x\in\re$ such that $v(x)\neq 0$. 
We refer to Theorem~\ref{thm:glaeser-k} below for a slightly more 
general result. Example~\ref{ex:glaeser-v'} and 
Example~\ref{ex:no-glaeser} below show the optimality of the 
assumptions. It is quite interesting that for $k\geq 2$ we have to 
impose that also $v'(x)$ does not change its sign, but we do not need 
similar assumptions on higher order derivatives (even if such 
assumptions would strongly simplify proofs).

The second question we address in this paper is to obtain similar
estimates for the first derivative of roots of a function defined in
some interval $(a,b)$.  So we take a function $g:(a,b)\to\re$, an
integer $k\geq 1$, and any continuous function $f:(a,b)\to\re$ such
that 
\begin{equation}
	|f(x)|^{k+1}=|g(x)|\quad\quad\forall x\in(a,b).
	\label{eqn:f-g}
\end{equation}

Trivial examples suggest that without assumptions on the sign of $g$
we cannot expect more than some sort of absolute continuity of $f$.

This problem was considered for the first time in Lemma~1
of~\cite{CJS}, where the authors proved that $f$ is absolutely
continuous in $(a,b)$ provided that $g\in C^{k+1}((a,b))$ is a
\emph{nonnegative} function.  This is a quite strong requirement
which, for example, forces $g'$ to vanish whenever $g$ vanishes.  More
recently, S.~Tarama~\cite{tarama} extended the technique introduced
in~\cite{CJS} and proved the same conclusion without assuming that $g$
is nonnegative.  We state Tarama's result in
section~\ref{sec:statements} as Theorem~\ref{thmbibl:tarama}.

The absolute continuity of $f$ is equivalent to saying that $f'\in
L^{1}((a,b))$.  In this paper we prove (Theorem~\ref{thm:main}) a
better summability of $f'$, namely that $f'\in L^{p}((a,b))$ for every
$p<1+1/k$.  Example~\ref{ex:p} and Example~\ref{ex:holder} below show
the optimality of our result, both with respect to $p$, and with
respect to the regularity of $g$.  Although the result is one
dimensional, it can be easily extended to any space dimension through
a straightforward sectioning argument (see
Theorem~\ref{thm:main-ndim}).

A higher summability of $f'$ was obtained by F.\ Colombini and N.\
Lerner in~\cite{colombini-lerner}, but once again in the case where
$g$ is nonnegative.  They proved that for every nonnegative $g\in
C^{k+1}((a,b))$ one has that $f'\in L^{p}((a,b))$ (locally) for every
$p<1+2/(k-1)$, and the same in any space dimension.  Of course all
solutions of~(\ref{eqn:f-g}) are also solutions of
$|f(x)|^{2k+2}=|g(x)|^{2}$, and thus their result implies that $f'\in
L^{p}((a,b))$ for every $p<1+1/k$ provided that $g\in
C^{2k+2}((a,b))$, without sign restrictions on $g$.  In other words,
one obtains our summability of $f'$ but with twice stronger regularity
assumptions on $g$.  This is observed in a remark at the end of section~4
of~\cite{colombini-lerner}.  The authors conclude the same remark by
pointing out that it is plausible that the only assumption $g\in C^{k+1}$ is
needed in order to obtain the same summability, but the proof of this fact
would require a nontrivial modification of their arguments.

In this paper we prove this conjecture, and indeed we follow a
completely different path with respect to the previous literature.
The technique used in~\cite{CJS} and~\cite{tarama} is oriented to
obtaining estimates of the total variation, namely \emph{integral}
estimates.  In our approach the higher summability of $f'$ follows
from suitable \emph{pointwise} estimates on $f'(x)$.  To this end, we
first divide $(a,b)$ into intervals $(a_{i},b_{i})$ where $g(x)\cdot
g'(x)\neq 0$.  Then in each interval we exploit a natural
generalization of Glaeser inequalities~(\ref{glaeser-k}) to functions
defined in bounded intervals (see Proposition~\ref{prop:ab}).  Thanks
to this generalization, we have an estimate on $f'(x)$, multiplied by
a standard cutoff function, in terms of the derivatives of $g$ up to
order $(k+1)$, and of the length $(b_{i}-a_{i})$ of the interval
itself.  These estimates yield an upper bound on the measure of the
sublevels of $f'$ in $(a_{i},b_{i})$.  Summing over all subintervals
we obtain the conclusion.

We point out that this approach works because the set of points where
$g(x)\cdot g'(x)=0$ does not contribute to the integral of $f'(x)$.
We stress that this would be false for the set of points where
$g(x)\cdot g'(x)\cdot g''(x)=0$, and this is the reason why it is
important to obtain higher order Glaeser inequalities such as
(\ref{glaeser-k}) with assumptions on $v$ and $v'$ only.

We conclude by mentioning some related results.  Regularity of roots
has long been studied, both because the problem is interesting in
itself, and because of several applications, ranging from algebraic
geometry to partial differential equations, for which we refer to the
quoted references.

A first research line concerns square roots, namely the case where
$f^{2}=g$ (see~\cite{BCP-sns} and the references quoted therein).  The
leitmotif of these papers is the search of conditions on $g$ under
which one can \emph{choose} a regular enough root $f$.  The usual
requirement is the existence of sufficiently many derivatives of $g$
and their vanishing at the zeroes of $g$.

Our problem is different for two reasons.  Firstly we investigate the
regularity of a \emph{given} root, which we cannot choose.  Secondly,
we have no assumptions on the sign of $g$, since we are actually
solving $|f|^{k+1}=|g|$.

Roots of functions are just the first step toward roots of
\emph{polynomials} with coefficients depending on a parameter.
Regularity results for roots of nonnegative functions have their
counterpart in analogous results for roots of hyperbolic polynomials,
namely polynomials whose roots are all real.  Classical conclusions
are that one can \emph{choose} roots of class $C^{1}$ or twice
differentiable provided that coefficients are smooth enough
(see~\cite{KLM1,bronstein,COP1,KLM2,mandai,tarama-2,waka}).

The situation considered in this paper corresponds to non-hyperbolic
polynomials, in which case one investigates the regularity of a given
continuous root.  The problem is still quite open.  T.\
Kato~\cite{kato} proved that for a monic polynomial of degree $n$ with
continuous coefficients one can always choose $n$ continuous functions
representing the roots of the polynomial.  S.\
Spagnolo~\cite{spagnolo-ferrara} proved that continuous roots of monic
polynomials of degree~2 and 3 are locally absolutely continuous
provided that coefficients are of class $C^{5}$ and $C^{25}$,
respectively.  This result is based on the explicit formulae for
solutions, and on the nonoptimal Lemma~1 of~\cite{CJS}, and for this
reason it is probably nonoptimal.

Spagnolo's result was in some sense extended by A.\
Rainer~\cite{rainer}, who proved the absolute continuity of roots of
monic polynomials of arbitrary degree with coefficients of class
$C^{\infty}$, but with the assumption that one can continuously
arrange roots in such a way that no two of them meet with an infinite
order of flatness.  In the case of the polynomial $P(y):=y^{k}-g(x)$,
namely in the case of $k$-th roots, this implies that
$g$ has a finite number of zeroes.  This simplifying assumption
greatly reduces the difficulty of the problem, but the estimates
obtained in this way do blow up when the number of zeroes goes to
infinity.  In this paper we do not impose this simplifying assumption,
and indeed the potential accumulation of zeroes of $g$ and its
derivatives is the main difficulty we have to face in our proofs.

We hope that our new approach based on pointwise estimates could be
helpful in order to understand in full generality the case of non-hyperbolic
polynomials.

This paper is organized as follows. In section~\ref{sec:statements} 
we fix notations, and we state our main results. In 
section~\ref{sec:proofs} we prove them. In section~\ref{sec:ex} we 
present the examples showing the optimality of our results, and the 
open problem concerning the regularity of roots of non-hyperbolic 
polynomials.

\setcounter{equation}{0}
\section{Notations and statements}\label{sec:statements}

Let $\Omega\subseteq\re$ be an open set, and let $\alpha\in(0,1]$.  A
function $f:\Omega\to\re$ is called \aholder\ continuous if there
exists a constant $H$ such that
\begin{equation}
	|f(y)-f(x)|\leq H|y-x|^{\alpha}
	\quad\quad
	\forall (x,y)\in\Omega^{2}.
	\label{defn:holder}
\end{equation}

In this case
$$\hold\left(f,\Omega\right):=\sup\left\{
\frac{|f(y)-f(x)|}{|y-x|^{\alpha}}:(x,y)\in\Omega^{2},\ x\neq 
y\right\}$$
denotes the \aholder\ constant of $f$ in $\Omega$, namely the 
smallest $H$ for which (\ref{defn:holder}) holds true.
Note that this includes the case of Lipschitz continuous functions
($\alpha=1$).

If $k$ is a positive integer, $C^{k,\alpha}(\Omega)$ denotes the set 
of functions $g\in C^{k}(\Omega)$ with derivatives up to order $k$ 
which are bounded in $\Omega$, and whose $k$-th derivative $g^{(k)}$ 
is \aholder\ continuous in $\Omega$. 

Regularity of roots of elements of $C^{k,\alpha}$ was investigated
in~\cite{tarama}.  The precise result is the following.

\begin{thmbibl}\label{thmbibl:tarama}
	Let $k$ be a positive integer, let $\alpha\in(0,1]$, let
	$(a,b)\subseteq\re$ be an interval, and let $f:(a,b)\to\re$ and 
	$g:(a,b)\to\re$ be two
	functions.  
	
	Let us assume that $f$ is continuous in $(a,b)$, and 
	that $g\in C^{k,\alpha}((a,b))$ satisfies
	\begin{equation}
		|f(x)|^{k+\alpha}=|g(x)|
		\quad\quad
		\forall x\in (a,b).
		\label{hp:fg}
	\end{equation}
	
	Then we have that $f'\in L^{1}((a,b))$, and there exists a 
	constant $C(k)$ such that
	$$\|f'\|_{L^{1}((a,b))}\leq
	C(k)\left\{\hold\left(g^{(k)},(a,b)\right)
	(b-a)^{k+\alpha}+
	\sum_{i=1}^{k}\|g^{(i)}\|_{L^{\infty}((a,b))}
	(b-a)^{i}\right\}^{1/(k+\alpha)}.$$
\end{thmbibl}

We are now ready to state our main contributions.  The first result of
this paper is the following generalization of Glaeser inequality.

\begin{thm}[Higher order Glaeser inequalities]\label{thm:glaeser-k}
	Let $k$ be a positive integer, let $\alpha\in(0,1]$, and let $v\in
	C^{k}(\re)$ be a function such that $v(x)$ and $v'(x)$ do not
	change their sign in $\re$ (namely either $v(x)\geq 0$ or
	$v(x)\leq 0$ for every $x\in\re$, and either $v'(x)\geq 0$ or
	$v'(x)\leq 0$ for every $x\in\re$).  
	
	Let us assume that the $k$-th derivative $v^{(k)}(x)$ is \aholder\
	continuous in $\re$.
	
	Then there exists a constant $C(k)$, depending only upon $k$, 
	such that
	\begin{equation}
		|v'(x)|^{k+\alpha}\leq C(k)\cdot|v(x)|^{k+\alpha-1}\cdot
		\hold\left(v^{(k)},\re\right)
		\quad\quad
		\forall x\in\re.
		\label{th:glaeser-k}
	\end{equation}
\end{thm}

A careful inspection of the proof of Theorem~\ref{thm:glaeser-k}
reveals that in the case $k=1$ we do not need assumptions on the sign
of $v'(x)$, as in the standard Glaeser inequality (\ref{glaeser}).

Theorem~\ref{thm:glaeser-k} above is optimal. Indeed
\begin{itemize}
	\item in Example~\ref{ex:glaeser-v'} we construct a positive
	function $v\in C^{\infty}(\re)$, whose derivatives are all 
	bounded, and such that the ratio
	\begin{equation}
		\frac{|v'(x)|^{k+\alpha}}{[v(x)]^{k+\alpha-1}}
		\label{ratio}
	\end{equation}
	is unbounded for every $k\geq 2$ and every $\alpha>0$ (of course
	$v'(x)$ is a sign changing function),

	\item in Example~\ref{ex:no-glaeser} we construct a function $v\in
	C^{\infty}(\re)$ such that the product $v(x)\cdot v'(x)$ never
	vanishes, $v\in C^{k,\beta}(\re)$ for every $\beta<\alpha$, and
	the ratio (\ref{ratio}) is unbounded.
\end{itemize}

In the second result we consider functions defined on an interval
$(a,b)$, and we improve Theorem~\ref{thmbibl:tarama} by showing a
better summability of $f'$.  As in~\cite{colombini-lerner} we express
this better summability in terms of weak $L^{p}$ spaces, whose
definition we briefly recall.

If $\Omega\subseteq\re$ is a bounded open set, and 
$\psi:\Omega\to\re$ is a measurable function, then for every $p\geq 1$ 
one can set
$$\|\psi\|_{p,w,\Omega}:=\sup_{M\geq 0}\left\{M\cdot\left[\strut
\mathrm{meas}\{x\in\Omega:|\psi(x)|>M\}\right]^{1/p}\right\},$$
where ``$\mathrm{meas}$'' denotes the Lebesgue measure. The weak $L^{p}$ 
space $L^{p}_{w}(\Omega)$ is the set of all functions for which this 
quantity (which is not a norm, since the triangular inequality fails 
to hold) is finite. In Lemma~\ref{lemma:weak-lp} we state all the 
properties of weak $L^{p}$ spaces needed in this paper. Here we just 
recall that for every $q<p$ we have that
$L^{p}(\Omega)\subset L^{p}_{w}(\Omega)\subset L^{q}(\Omega)$,
with strict inclusions.

We are now ready to state our main result.

\begin{thm}[Optimal regularity of roots]\label{thm:main}
	Let $k$ be a positive integer, let $\alpha\in(0,1]$, let
	$(a,b)\subseteq\re$ be an interval, and let $f:(a,b)\to\re$ be a
	function.  
	
	Let us assume that $f$ is continuous in $(a,b)$, and
	there exists $g\in C^{k,\alpha}((a,b))$ such that (\ref{hp:fg})
	holds true. Let $p=p(k,\alpha)$ be such that
	\begin{equation}
		\frac{1}{p}+\frac{1}{k+\alpha}=1.
		\label{defn:p*}
	\end{equation}
	
	Then we have that $f'\in L^{p}_{w}((a,b))$, and there exists a
	constant $C(k)$, depending only upon $k$, such that
	\begin{equation}
		\|f'\|_{p,w,(a,b)} \leq C(k)
		\max\left\{
		\left[\hold\left(g^{(k)},(a,b)\right)\right]^{1/(k+\alpha)}
		\cdot(b-a)^{1/p},
		\|g'\|^{1/(k+\alpha)}_{L^{\infty}((a,b))}\right\}.
		\label{th:main}
	\end{equation}
	
	In particular, we have that $f$ lies in the Sobolev space
	$W^{1,q}((a,b))$ for every $q\in[1,p)$.
\end{thm}

We point out that in the right-hand side of (\ref{th:main}) we have
the maximum between two terms which are homogeneous of degree $1-1/p$
with respect to dilatations (namely when replacing $g(x)$ with
$g(\lambda x)$).  Also the left-hand side of~(\ref{th:main}) is
homogeneous with the same degree.  Moreover Theorem~\ref{thm:main}
above is optimal, in the sense that
\begin{itemize}
	\item in Example~\ref{ex:p} below we construct $f$ and $g$ such
	that $g\in C^{k,\alpha}((a,b))\cap C^{\infty}((a,b))$, but
	$f'\not\in L^{p}((a,b))$ for $p$ given by (\ref{defn:p*}),

	\item in Example~\ref{ex:holder} below we construct $f$ and $g$
	such that $g\in C^{\infty}((a,b))$ and $g\in C^{k,\beta}((a,b))$
	for every $\beta<\alpha$, but $f$ is not a bounded variation
	function in $(a,b)$, hence also $f'\not\in L^{1}((a,b))$.
\end{itemize}

Example~\ref{ex:holder} shows also that there is some sort of gap
phenomenon.  If $g\in C^{k,\beta}$ for all $\beta<\alpha$, then
it is possible that its $(k+\alpha)$-th root $f$ is not even a bounded
variation function.  As soon as $g\in C^{k,\alpha}$ we have
that $f\in W^{1,q}$ for every $q\in[1,p)$.  So we
have an immediate jump in the regularity of the root without passing
through intermediate values of the Sobolev exponent.

We conclude by extending Theorem~\ref{thm:main} to higher dimension.
The extension is straightforward because the spaces
$C^{k,\alpha}(\Omega)$ and $L^{p}_{w}(\Omega)$ can be easily defined
for every open set $\Omega\subseteq\re^{n}$, and the Sobolev regularity
of a function of $n$ real variables depends on the Sobolev regularity
of its one-dimensional sections.

We obtain the following result.

\begin{thm}[Higher dimensional case]\label{thm:main-ndim}
	Let $n$ and $k$ be two positive integers, let $\alpha\in(0,1]$,
	let $\Omega\subseteq\re^{n}$ be an open set, and let
	$f:\Omega\to\re$ be a function.  
	
	Let us assume that $f$ is continuous in $\Omega$, and there exists
	$g\in C^{k,\alpha}(\Omega)$ such that (\ref{hp:fg}) holds true for
	every $x\in\Omega$. Let $p$ be defined by (\ref{defn:p*}), and 
	let $\hold\left(g^{(k)},\Omega\right)$ be the maximum of the 
	H\"{o}lder constants in $\Omega$ of all partial derivatives of 
	$g$ with order $k$.  
	
	Then $\nabla f\in L^{p}_{w}(\Omega')$ for every
	$\Omega'\subset\subset\Omega$, and there exist constants
	$C_{1}(n,k,\Omega,\Omega')$ and $C_{2}(n,k,\Omega,\Omega')$, both
	independent of $f$, such that
	$$\|\nabla f\|_{p,w,\Omega'}\leq
	\max\left\{
	C_{1}(n,k,\Omega,\Omega')
	\left[\hold\left(g^{(k)},\Omega\right)\right]^{1/(k+\alpha)},
	C_{2}(n,k,\Omega,\Omega')
	\|\nabla g\|_{L^{\infty}(\Omega)}^{1/(k+\alpha)}\right\}.$$
\end{thm}

A careful inspection of the proof reveals that $C_{1}$ and $C_{2}$
have homogeneity degrees with respect to dilatations equal to $-n/p$
and $-(n-1)/p$, respectively, which makes the two terms in the maximum
homogeneous with the same degree of the left-hand side.  It also
reveals that what we actually need to assume on $g$ is the regularity
of its one-dimensional sections obtained by freezing $(n-1)$ variables
in all possible ways.

\section{Proofs}\label{sec:proofs}

\subsection{Technical preliminaries}

\begin{lemma}[Properties of weak $L^{p}$ spaces]\label{lemma:weak-lp}
	For every bounded open set $\Omega\subseteq\re$, and every $p\geq 
	1$, we have the following conclusions.
	\begin{enumerate}
		\renewcommand{\labelenumi}{(\arabic{enumi})}
	
		\item  If $\psi_{1},\ldots,\psi_{n}$ are elements of 
		$L^{p}_{w}(\Omega)$, then
		$$\|\psi_{1}+\ldots+\psi_{n}\|_{p,w,\Omega}\leq n\left(\strut
		\|\psi_{1}\|_{p,w,\Omega}+\ldots+\|\psi_{n}\|_{p,w,\Omega}\right).$$
	
		\item  We have that 
		$$\|\lambda\psi\|_{p,w,\Omega}=|\lambda|\cdot\|\psi\|_{p,w,\Omega}
		\quad\quad
		\forall\lambda\in\re,\ \forall\psi\in L^{p}_{w}(\Omega).$$
		
		\item  For every $\psi_{1}$ and $\psi_{2}$ in 
		$L^{p}_{w}(\Omega)$ we have that
		$$\|\max\{\psi_{1},\psi_{2}\}\|_{p,w,\Omega}\leq
		2\max\{\|\psi_{1}\|_{p,w,\Omega},\|\psi_{2}\|_{p,w,\Omega}\}.$$
	
		\item Let $\Omega:=(a,b)$ be an interval, and let
		$\psi_{1}(x):=1$ and
		$\psi_{2}(x):=\left[(x-a)(b-x)\right]^{-1/p}$ for every
		$x\in(a,b)$.  Then we have that
		$$\|\psi_{1}\|_{p,w,(a,b)}=(b-a)^{1/p}, \hspace{4em}
		\|\psi_{2}\|_{p,w,(a,b)}\leq
		\left(\frac{4}{b-a}\right)^{1/p}.$$
	
		\item Let $\{\Omega_{J}\}_{J\in\mathcal{C}}$ be a finite or
		countable family of open sets whose union is $\Omega$.  Then
		we have that 
		$$\|\psi\|_{p,w,\Omega}^{p}\leq \sum_{J\in\mathcal{C}}
		\|\psi\|_{p,w,\Omega_{J}}^{p} \quad\quad \forall\psi\in
		L^{p}_{w}(\Omega).$$
	\end{enumerate}
\end{lemma}

\paragraph{\textmd{\emph{Proof}}}

All these properties are simple consequences of the definition of 
weak $L^{p}$ spaces. In order to estimate $\|\psi_{2}\|_{p,w,(a,b)}$ 
in statement~(4) it is enough to remark that
$$\hspace{-2em}
\mathrm{meas}\{x\in(a,b):|\psi_{2}(x)|>M\}=
\mathrm{meas}\{x\in(a,b):(x-a)(b-x)<M^{-p}\}$$
$$\hspace{2em}
\leq 2\,\mathrm{meas}
\left\{x\in\left(a,\frac{a+b}{2}\right):(x-a)\left(\frac{b-a}{2}\right)
<\frac{1}{M^{p}}\right\}\leq
\frac{4}{b-a}\cdot\frac{1}{M^{p}}.\qed$$

\begin{lemma}\label{lemma:AC}
	Let $(a,b)\subseteq\re$ be an interval, let $f:(a,b)\to\re$ be a
	continuous function, and let $\Omega_{0}:=\{x\in(a,b):f(x)\neq
	0\}$.  
	
	Let us assume that $f\in C^{1}(\Omega_{0})$, and that $f'\in 
	L^{p}_{w}(\Omega_{0})$ for some $p>1$.
	
	Then we have that the distributional derivative of $f$ in $(a,b)$
	is a measurable function $f'\in L^{p}_{w}((a,b))$, and
	\begin{equation}
		\|f'\|_{p,w,(a,b)}=\|f'\|_{p,w,\Omega_{0}}.
		\label{th:lemma-ac}
	\end{equation}
\end{lemma}

\paragraph{\textmd{\emph{Proof}}}

Let us set
$$\psi(x):=\left\{
\begin{array}{ll}
	f'(x) & \mbox{if }x\in\Omega_{0},  \\
	0 & \mbox{if }x\in(a,b)\setminus\Omega_{0}.
\end{array}
\right.$$

Clearly we have that $\psi\in L^{p}_{w}((a,b))$.  We claim that $\psi$
is the distributional derivative of $f$ in
$(a,b)$.  Indeed let us take any compactly supported test function
$\phi\in C^{\infty}((a,b))$.  Let $\mathcal{C}$ denote the (finite or
countable) set of connected components of $\Omega_{0}$.  Then in each
$J\in\mathcal{C}$ we can use the standard integration-by-parts formula
for smooth functions, and obtain that
$$\hspace{-4em}
\int_{a}^{b}f(x)\phi'(x)\,dx=\int_{\Omega_{0}}f(x)\phi'(x)\,dx=
\sum_{J\in\mathcal{C}}\int_{J}f(x)\phi'(x)\,dx$$
$$\hspace{4em}
=-\sum_{J\in\mathcal{C}}\int_{J}f'(x)\phi(x)\,dx=
-\int_{\Omega_{0}}f'(x)\phi(x)\,dx=
-\int_{a}^{b}\psi(x)\phi(x)\,dx.$$

Now that we know that $\psi$ is the distributional derivative of $f$ 
it is easy to conclude that
$$\|f'\|_{p,w,(a,b)}=\|\psi\|_{p,w,(a,b)}=\|\psi\|_{p,w,\Omega_{0}}=
\|f'\|_{p,w,\Omega_{0}},$$ 
which completes the proof.\qed

\subsection{Key estimates}

The following surprising property of polynomials is the key tool in 
the proof of our higher order Glaeser inequalities.  We stress that the
``magic constants'' $\delta_{k}$ and $\ep_{k}$ do not depend on
$P(x)$.

\begin{lemma}[Magic constants]\label{lemma:main-pol}
	For every positive integer $k$ there exist positive constants 
	$\delta_{k}$ and $\ep_{k}$, depending only upon $k$, with the following 
	properties.
	\begin{enumerate}
		\renewcommand{\labelenumi}{(\arabic{enumi})}
		\item For every polynomial $P(x)$, with degree less than or
		equal to $k$, such that $P(0)=1$ and $P'(0)=-1$, there exists
		$x\in[0,\delta_{k}]$ such that either $P(x)\leq -1$ or
		$P'(x)\geq 1$.
	
		\item  For every real number $A$, every $\alpha\in(0,1]$, and
		every polynomial $P(x)$, with degree less than or equal to
		$k$, such that $P(0)=1$ and $P'(0)=-1$, the following
		implication
		\begin{equation}
			\fbox{$\begin{array}{ll}
				P(x)+Ax^{k+\alpha}\geq 0 & \forall x\in[0,\delta_{k}]  \\
				\noalign{\vspace{1ex}}
				P'(x)-Akx^{k+\alpha-1}\leq 0\quad & \forall x\in[0,\delta_{k}]
			\end{array}$}
			\quad\Longrightarrow\quad
			\fbox{$A\geq\ep_{k}$}
			\label{th:impl-main}
		\end{equation}
		holds true.
	\end{enumerate}
	 
\end{lemma}

\paragraph{\textmd{\emph{Proof}}}

Let us assume that statement~(1) is false for some $k$. Then there 
exists a sequence of polynomials $P_{n}(x)$, with degree less than or 
equal to $k$, such that $P_{n}(0)=1$, $P_{n}'(0)=-1$, and
\begin{equation}
	P_{n}(x)\geq -1
	\quad\mbox{and}\quad
	P'_{n}(x)\leq 1
	\quad\quad
	\forall x\in[0,n].
	\label{est:pn}
\end{equation}

These assumptions easily imply that
$$-1\leq P_{n}(x)\leq P_{n}(0)+x= 1+x\leq 2$$
for every $x\in[0,1]$. In particular, the sequence $\{P_{n}(x)\}$ is 
bounded for at least $(k+1)$ distinct values of $x$, which we denote 
by $x_{1}$, \ldots, $x_{k}$, $x_{k+1}$.

As a consequence, there exists a polynomial $P_{\infty}(x)$, with
degree at most $k$, such that $P_{n}(x)\to P_{\infty}(x)$ up to
subsequences (not relabeled), in the sense that we have both
convergence of coefficients, and pointwise convergence.  This is a
well known fact which follows, for example, from the simple remark
that coefficients of a polynomial $P(x)$ with degree less than or
equal to $k$ depend in a linear way (through a Vandermonde matrix)
from the values of $P(x)$ in the $(k+1)$ points $x_{1}$, \ldots,
$x_{k}$, $x_{k+1}$.  Alternatively, it follows from the fact that the
maximum of $|P(x_{i})|$ for $i=1,\ldots,k+1$, and the maximum of the
absolute value of the coefficients of $P(x)$ are two norms on the
space of all polynomials with degree less than or equal to $k$, and
they are equivalent because the space is finite dimensional.

In any case, convergence of coefficients implies that
$P_{\infty}(0)=1$ and $P_{\infty}'(0)=-1$.  Pointwise convergence
implies that we can pass to the limit in (\ref{est:pn}) and obtain
that
\begin{equation}
	P_{\infty}(x)\geq -1
	\quad\mbox{and}\quad
	P_{\infty}'(x)\leq 1
	\quad\quad
	\forall x\geq 0.
	\label{est:p-infty}
\end{equation}

If the degree of $P_{\infty}(x)$ is 1, then $P_{\infty}(x)=1-x$, which
does not satisfy (\ref{est:p-infty}).  If the degree of
$P_{\infty}(x)$ is at least~2, then as $x\to +\infty$ we have that
both $P_{\infty}(x)$ and $P_{\infty}'(x)$ tend to $+\infty$ or
$-\infty$ according to the sign of the leading coefficient.  This
contradicts (\ref{est:p-infty}) also in this case, and completes the 
proof of statement~(1).

Now we claim that the conclusion of statement~(2) is true with
$$\ep_{k}:=\min\left\{
\frac{1}{\delta_{k}^{k+1}},
\frac{1}{k\delta_{k}^{k}}\right\}.$$

Indeed from statement~(1) we know that there exists $x\in[0,\delta_{k}]$ 
such that either $P(x)\leq -1$ or $P'(x)\geq 1$. In the first case 
the first inequality in the left-hand side of (\ref{th:impl-main}) 
gives that (clearly $\delta_{k}\geq 1$ for every $k\geq 1$)
$$A\geq\frac{-P(x)}{x^{k+\alpha}}\geq\frac{1}{x^{k+\alpha}}\geq
\frac{1}{\delta_{k}^{k+\alpha}}\geq\frac{1}{\delta_{k}^{k+1}}.$$

In the second case the second inequality in the left-hand side of
(\ref{th:impl-main}) gives that
$$A\geq\frac{P'(x)}{kx^{k+\alpha-1}}\geq\frac{1}{kx^{k+\alpha-1}}\geq
\frac{1}{k\delta_{k}^{k+\alpha-1}}\geq\frac{1}{k\delta_{k}^{k}}.$$

In both cases implication (\ref{th:impl-main}) holds true.\qed
\medskip

Next result is the main step toward higher order Glaeser inequalities,
both in $\re$ and in bounded intervals.  The proof exploits the
``magic constants'' of Lemma~\ref{lemma:main-pol} and Taylor's
expansion.

\begin{prop}[Higher order Glaeser inequality in an interval]\label{prop:pre-glaeser}
	\hskip 0em plus 1em
	Let $k$ be a positive integer, and let $\alpha\in(0,1]$,
	$x_{0}\in\re$, and $d>0$ be real numbers.  
	
	Let $v\in
	C^{k}((x_{0}-d,x_{0}+d))$ be a function such that $v(x)$ and
	$v'(x)$ do not change their sign, namely either $v(x)\geq
	0$ or $v(x)\leq 0$ for every $x\in(x_{0}-d,x_{0}+d)$, and the same
	for $v'(x)$.
	Let us assume that the $k$-th derivative $v^{(k)}(x)$ is \aholder\
	continuous in $(x_{0}-d,x_{0}+d)$.
	
	Then there exists a constant $C(k)$, depending only upon $k$, 
	such that
	\begin{equation}
		|v'(x_{0})|^{k+\alpha}\leq
		C(k)\cdot|v(x_{0})|^{k+\alpha-1}\cdot \max\left\{
		\hold\left(v^{(k)},(x_{0}-d,x_{0}+d)\right),
		\frac{|v'(x_{0})|}{d^{k+\alpha-1}}\right\}.
		\label{th:pre-glaeser}
	\end{equation}
\end{prop}

\paragraph{\textmd{\emph{Proof}}}

Up to symmetries we can assume, without loss of generality, that
$v(x)\geq 0$ and $v'(x)\leq 0$ for every $x\in(x_{0}-d,x_{0}+d)$.  If
$v'(x_{0})=0$, then there is nothing to prove.  If $v(x_{0})=0$, then
it is easy to see that also $v'(x_{0})=0$, and once again there is
nothing to prove.  So from now on we assume that $v(x_{0})>0$ and
$v'(x_{0})<0$.

Let $\delta_{k}$ be the constant of Lemma~\ref{lemma:main-pol}. We 
distinguish two cases.

\subparagraph{\textmd{\textit{Case 1}}}

Let us assume that
$$\frac{|v'(x_{0})|}{v(x_{0})}\,d\leq\delta_{k}.$$

This is equivalent to saying that 
$|v'(x_{0})|\leq\delta_{k}d^{-1}v(x_{0})$, and hence
$$|v'(x_{0})|^{k+\alpha}=|v'(x_{0})|\cdot|v'(x_{0})|^{k+\alpha-1}\leq
|v'(x_{0})|\cdot
\left(\frac{\delta_{k}}{d}\right)^{k+\alpha-1}[v(x_{0})]^{k+\alpha-1},$$
which implies (\ref{th:pre-glaeser}) in this first case, provided 
that $C(k)\geq\delta_{k}^{k}$.

\subparagraph{\textmd{\textit{Case 2}}}

Let us assume now that
\begin{equation}
	\frac{|v'(x_{0})|}{v(x_{0})}\,d>\delta_{k}.
	\label{hp:case2}
\end{equation}

Let us set for simplicity
$H:=\hold\left(v^{(k)},(x_{0}-d,x_{0}+d)\right)$, and let us write
Taylor's expansion of $v(x)$ of order $k$.  For every
$\delta\in[0,d)$ we obtain that
$$v(x_{0}+\delta)=\sum_{j=0}^{k-1}\frac{v^{(j)}(x_{0})}{j!}\,\delta^{j}+
\frac{v^{(k)}(x_{0}+\xi)}{k!}\,\delta^{k}$$
for some $\xi\in(0,\delta)$. Since $v^{(k)}$ is \aholder\ continuous 
we have that
$$|v^{(k)}(x_{0}+\xi)-v^{(k)}(x_{0})|\leq H\delta^{\alpha}.$$

Therefore our assumption that $v(x)\geq 0$ for every
$x\in(x_{0}-d,x_{0}+d)$ implies that
$$0\leq v(x_{0}+\delta)\leq
\sum_{j=0}^{k}\frac{v^{(j)}(x_{0})}{j!}\,\delta^{j}+
\frac{H}{k!}\,\delta^{k+\alpha}
\quad\quad
\forall\delta\in[0,d).$$

In an analogous way we obtain that
$$0\geq v'(x_{0}+\delta)\geq
\sum_{j=0}^{k-1}\frac{v^{(j+1)}(x_{0})}{j!}\,\delta^{j}-
\frac{H}{(k-1)!}\,\delta^{k+\alpha-1}
\quad\quad
\forall\delta\in[0,d).$$

Let us write both inequalities with 
$\delta:=x\cdot v(x_{0})\cdot|v'(x_{0})|^{-1}$. We obtain that the 
two inequalities
$$0\leq|v(x_{0})|\left(\sum_{j=0}^{k}
\frac{v^{(j)}(x_{0})}{j!}\cdot
\frac{[v(x_{0})]^{j-1}}{|v'(x_{0})|^{j}}\cdot x^{j}+
\frac{H}{k!}\cdot
\frac{[v(x_{0})]^{k+\alpha-1}}{|v'(x_{0})|^{k+\alpha}}\cdot
x^{k+\alpha}\right),$$
$$0\geq|v'(x_{0})|\left(\sum_{j=0}^{k-1}
\frac{v^{(j+1)}(x_{0})}{j!}\cdot
\frac{[v(x_{0})]^{j}}{|v'(x_{0})|^{j+1}}\cdot x^{j}-
\frac{H}{(k-1)!}\cdot
\frac{[v(x_{0})]^{k+\alpha-1}}{|v'(x_{0})|^{k+\alpha}}\cdot
x^{k+\alpha-1}\right)$$
hold true for every 
$$0\leq x<\frac{|v'(x_{0})|}{v(x_{0})}\,d.$$ 

Due to (\ref{hp:case2}), the upper bound is larger than $\delta_{k}$, 
which proves that the two inequalities hold true at least for every 
$x\in[0,\delta_{k}]$. Thus if we set
\begin{equation}
	P(x):=\sum_{j=0}^{k}\frac{v^{(j)}(x_{0})}{j!}\cdot
	\frac{[v(x_{0})]^{j-1}}{|v'(x_{0})|^{j}}\cdot x^{j},
	\hspace{3em}
	A:=\frac{H}{k!}\cdot
	\frac{[v(x_{0})]^{k+\alpha-1}}{|v'(x_{0})|^{k+\alpha}},
	\label{defn:aj-A}
\end{equation}
we have that $P(0)=1$, $P'(0)=-1$, and the assumptions in the
left-hand side of (\ref{th:impl-main}) are satisfied.

Therefore from (\ref{th:impl-main}) it
follows that $A\geq\ep_{k}$, which implies
(\ref{th:pre-glaeser}) also in this second case, provided that 
$C(k)\geq(\ep_{k}k!)^{-1}$.\qed

\subsection{Proof of Theorem~\ref{thm:glaeser-k}}

Let us fix any point $x_{0}\in\re$.  Then the assumptions of
Proposition~\ref{prop:pre-glaeser} are satisfied for every $d>0$.  The
conclusion follows by remarking that
$$\hold\left(v^{(k)},(x_{0}-d,x_{0}+d)\right)\leq
\hold\left(v^{(k)},\re\right),$$
and finally letting $d\to+\infty$.\qed

\subsection{Absolute continuity of roots}

The aim of this section is to prove Theorem~\ref{thm:main} and 
Theorem~\ref{thm:main-ndim}. The proof is based on some \emph{pointwise} 
estimates on the derivative, which are the content of next two results.

\begin{prop}\label{prop:ab}
	Let $k$ be a positive integer, let $\alpha\in(0,1]$, let
	$(a,b)\subseteq\re$ be an interval, let $f:(a,b)\to\re$ be a
	continuous function, and let $g\in C^{k,\alpha}((a,b))$.  
	
	Let us assume that $f$ and $g$ satisfy (\ref{hp:fg}), and that
	$g(x)\cdot g'(x)\neq 0$ for every $x\in (a,b)$.  Let us set for
	simplicity $$H:=\hold\left(g^{(k)},(a,b)\right), \hspace{3em}
	\|g'\|_{\infty}:=\|g'\|_{L^{\infty}((a,b))}.$$
	
	Then we have the following conclusions.
	\begin{enumerate}
		\renewcommand{\labelenumi}{(\arabic{enumi})}
		\item  There exists a constant $C_{1}(k)$, depending only upon 
		$k$, such that for every $x\in(a,b)$ we have that
		\begin{equation}
			|f'(x)|^{k+\alpha}\leq C_{1}(k)\max\left\{
			H,
			\|g'\|_{\infty}\cdot\left[\frac{b-a}
			{(x-a)(b-x)}\right]^{k+\alpha-1}\right\}.
			\label{th:ab}
		\end{equation}

		\item Let $p$ be defined by (\ref{defn:p*}). 
		Then there exists a constant $C_{2}(k)$, depending only upon 
		$k$, such that
		\begin{equation}
			\|f'\|_{p,w,(a,b)}\leq C_{2}(k)\max\left\{\strut
			H^{1/(k+\alpha)}\cdot
			(b-a)^{1/p},
			\|g'\|_{\infty}^{1/(k+\alpha)}\right\}.
			\label{th:ab-fp}
		\end{equation}

	\end{enumerate}
\end{prop}

\paragraph{\textmd{\emph{Proof}}}

Let us fix any point $x_{0}\in(a,b)$, and let
$d:=\min\{x_{0}-a,b-x_{0}\}$. We can apply
Proposition~\ref{prop:pre-glaeser} to the function $v(x):=g(x)$ in the
interval $(x_{0}-d,x_{0}+d)$.  We obtain that
$$|f'(x_{0})|^{k+\alpha}=
\frac{|g'(x_{0})|^{k+\alpha}}{|g(x_{0})|^{k+\alpha-1}}\leq
C(k)\max\left\{H,
\frac{|g'(x_{0})|}{d^{k+\alpha-1}}\right\}.$$

Since
$$\frac{1}{d^{k+\alpha-1}}\leq\left[\frac{b-a}
{(x_{0}-a)(b-x_{0})}\right]^{k+\alpha-1},$$
conclusion (\ref{th:ab}) easily follows.

Let us consider now statement~(2).
From (\ref{th:ab}) and definition~(\ref{defn:p*}) of $p$, we have that 
$|f'(x)|\leq\left[C_{1}(k)\right]^{1/(k+\alpha)}
\max\{\psi_{1}(x),\psi_{2}(x)\}$, with
$$\psi_{1}(x):=
H^{1/(k+\alpha)},
\hspace{4em}
\psi_{2}(x):=\|g'\|_{\infty}^{1/(k+\alpha)}\cdot
\frac{(b-a)^{1/p}}{\left[(x-a)(b-x)\right]^{1/p}}.$$

From statements (2) through (4) of Lemma~\ref{lemma:weak-lp} we easily 
obtain (\ref{th:ab-fp}).\qed

\begin{prop}\label{prop:cd}
	Let $k$ be a positive integer, let $\alpha\in(0,1]$, let
	$(a,b)\subseteq[c,d]\subseteq\re$ be an open and a closed
	interval, let $f:(a,b)\to\re$ be a continuous function, and let
	$g\in C^{k}([c,d])$.  Let us assume that
	\begin{enumerate}
		\renewcommand{\labelenumi}{(\roman{enumi})} 
		\item $f$ and $g$ satisfy (\ref{hp:fg}),
		
		\item $g(x)\cdot g'(x)\neq 0$ for every $x\in (a,b)$,
		
		\item the $k$-th derivative $g^{(k)}(x)$ is \aholder\ 
		continuous in $(c,d)$ (hence also in $[c,d]$),
		
		\item for every $h=1,\ldots,k$ there exists 
		$x_{h}\in[c,d]$ such that $g^{(h)}(x_{h})=0$.
	\end{enumerate}
	
	Let $p$ be defined by (\ref{defn:p*}), and let $C_{2}(k)$ the the
	constant in~(\ref{th:ab-fp}).  Then we have that
	\begin{equation}
		\|f'\|_{p,w,(a,b)}\leq C_{2}(k)
		\left[\hold\left(g^{(k)},(c,d)\right)\right]^{1/(k+\alpha)}
		\cdot(d-c)^{1/p}.
		\label{th:cd-fp}
	\end{equation}
\end{prop}

\paragraph{\textmd{\emph{Proof}}}

Let us set for simplicity $H:=\hold\left(g^{(k)},(c,d)\right)$.  Due
to assumptions (i), (ii), and (iii), we can apply
Proposition~\ref{prop:ab} and obtain estimate (\ref{th:ab-fp}).  Due to
assumptions (iii) and (iv), we can estimate in $[c,d]$ all derivatives
of $g$ up to order $k$.  For every $h=1,\ldots,k$ we obtain that
$$|g^{(h)}(x)|\leq H\cdot (d-c)^{k+\alpha-h}
\quad\quad
\forall x\in[c,d]$$
(formally one should argue by induction on $k-h$), and in particular
\begin{equation}
	\|g'\|_{L^{\infty}((a,b))}\leq 
	\|g'\|_{L^{\infty}((c,d))}\leq H\cdot(d-c)^{k+\alpha-1}.
	\label{est:cd}
\end{equation}

Plugging this inequality into (\ref{th:ab-fp}), and estimating 
$(b-a)$ with $(d-c)$, we easily obtain (\ref{th:cd-fp}).\qed

\subsubsection*{Proof of Theorem~\ref{thm:main}}

Let us consider the open set
$$\Omega_{0}:=\{x\in(a,b):f(x)\neq 0\}=\{x\in(a,b):g(x)\neq 0\}.$$

From (\ref{hp:fg}) it is easy to see that $f\in C^{1}(\Omega_{0})$ (or
better $f\in C^{k}(\Omega_{0})$), and
$$\Omega_{1}:=\{x\in\Omega_{0}:f'(x)\neq 0\}=
\{x\in\Omega_{0}:g'(x)\neq 0\}.$$

Let $\mathcal{C}$ denote the
(finite or countable) set of connected components of $\Omega_{1}$.
Each $J\in\mathcal{C}$ is an open interval of the form
$(a_{J},b_{J})$, so that
$$\Omega_{1}=\bigcup_{J\in \mathcal{C}}J=
\bigcup_{J\in \mathcal{C}}(a_{J},b_{J}).$$

Thanks to Lemma~\ref{lemma:AC} and statement~(5) of 
Lemma~\ref{lemma:weak-lp} we have that
\begin{equation}
	\|f'\|_{p,w,(a,b)}^{p}= \|f'\|_{p,w,\Omega_{0}}^{p}
	= \|f'\|_{p,w,\Omega_{1}}^{p}
	\leq \sum_{J\in \mathcal{C}} \|f'\|_{p,w,J}^{p}.
	\label{th:main-bis}
\end{equation}

It is therefore enough to show that the right-hand side of
(\ref{th:main-bis}) is estimated by the right-hand side of
(\ref{th:main}), raised to the power of $p$, for a suitable constant
$C(k)$ depending only upon $k$.  To this end, we divide the connected
components of $\Omega_{1}$ into two disjoint classes $\mathcal{C}_{0}$
and $\mathcal{C}_{1}$.

\subparagraph{\textmd{\emph{Partitioning connected components}}}

Let $\mathcal{C}_{0}\subseteq \mathcal{C}$ be the set of connected
components $J\in \mathcal{C}$ for which $(a_{J},b_{J})$ is such
that either $a_{J}=a$ or
\begin{equation}
	\left|[b_{J},b)\setminus\Omega_{1}\strut\right|\leq 2k,
	\label{defn:c0}
\end{equation}
where vertical bars denote the number of elements of a set, and let
$\mathcal{C}_{1}:=\mathcal{C}\setminus \mathcal{C}_{0}$.  If
$(a,b)\setminus\Omega_{1}$ is a finite set, then $\mathcal{C}$ is
finite and $\mathcal{C}_{0}$ contains the leftmost and the $(2k+1)$
rightmost connected components of $\Omega_{1}$, or
$\mathcal{C}_{0}=\mathcal{C}$ if $|\mathcal{C}|\leq 2k+1$.

If $(a,b)\setminus\Omega_{1}$ is not finite, then $\mathcal{C}_{0}$
contains the connected component of $\Omega_{1}$ with $a_{J}=a$ (if it
exists), and the rightmost $h$ connected components of $\Omega_{1}$,
where $h$ is the largest integer less than or equal to $(2k+1)$ such
that the rightmost $h$ connected components of $\Omega_{1}$ exist and
their right-hand endpoints are not accumulation points of
$(a,b)\setminus\Omega_{1}$.

Note that it may happen that $a$ and $b$ themselves are accumulation
points of $(a,b)\setminus\Omega_{1}$, and in this case
$\mathcal{C}_{0}$ is empty.  It is also possible that
$|\mathcal{C}|=2k+2$ but $|\mathcal{C}_{0}|<2k+2$, because it could happen
that $(a,b)\setminus\Omega_{1}$ contains an interval.

What is important is that in any case we have that
$|\mathcal{C}_{0}|\leq 2k+2$.

\subparagraph{\textmd{\emph{Absolute continuity in the union of
$\mathcal{C}_{0}$-components}}}

Let us apply Proposition~\ref{prop:ab} in a connected component $J\in
\mathcal{C}_{0}$.  We obtain that
\begin{eqnarray*}
	\|f'\|^{p}_{p,w,J} & \leq &
	[C_{2}(k)]^{p}\max\left\{
	\left[\hold\left(g^{(k)},J)\right)\right]^{p/(k+\alpha)}
	\cdot (b_{J}-a_{J}),
	\|g'\|_{L^{\infty}(J)}^{p/(k+\alpha)}\right\} \\
	\noalign{\vspace{1ex}}
	 & \leq & [C_{2}(k)]^{p}\max\left\{
	\left[\hold\left(g^{(k)},(a,b)\right)\right]^{p/(k+\alpha)}
	\cdot (b-a),
	\|g'\|_{L^{\infty}((a,b))}^{p/(k+\alpha)}\right\}.
\end{eqnarray*}

Summing over all $J\in \mathcal{C}_{0}$, and recalling that
$|\mathcal{C}_{0}|\leq 2k+2$, we deduce that
\begin{eqnarray}
	\lefteqn{\hspace{-3em}
	\sum_{J\in \mathcal{C}_{0}}\|f'\|^{p}_{p,w,J}  
	\ \leq\  
	(2k+2)\cdot [C_{2}(k)]^{p}\cdot}
	\nonumber \\
	\quad\quad && \cdot\max\left\{
	\left[\hold\left(g^{(k)},(a,b)\right)\right]^{p/(k+\alpha)} \cdot
	(b-a), \|g'\|_{L^{\infty}((a,b))}^{p/(k+\alpha)}\right\}.
	\label{est:c0}
\end{eqnarray}

\subparagraph{\textmd{\emph{Expanding $\mathcal{C}_{1}$-components}}}

We prove that to every $J\in \mathcal{C}_{1}$ we can associate a
closed interval $[c_{J},d_{J}]$ in such a way that the following three
conditions are satisfied.
\begin{enumerate}
	\renewcommand{\labelenumi}{(P\arabic{enumi})}
	\item  We have that 
	$(a_{J},b_{J})\subseteq[c_{J},d_{J}]\subseteq(a,b)$ for 
	every $J\in \mathcal{C}_{1}$.

	\item  Every $x\in (a,b)$ lies in the interior of at most 
	$(2k+1)$ such intervals, namely
	\begin{equation}
		\left|\strut\{J\in \mathcal{C}_{1}:x\in(c_{J},d_{J})\}\right|
		\leq 2k+1.
		\label{eqn:P2}
	\end{equation}

	\item For every $J\in \mathcal{C}_{1}$, and every $h=1,\ldots,k$,
	there exists $x_{J,h}\in[c_{J},d_{J}]$ such that
	$g^{(h)}(x_{J,h})=0$ (note that $g$ is always defined in the 
	closed interval $[c_{J},d_{J}]$).
\end{enumerate}

To this end, for every $J\in \mathcal{C}_{1}$ we set
$$D_{J}:=\left\{x\in[b_{J},b):
\left|[b_{J},x]\setminus\Omega_{1}\strut\right|\geq 2k+1\right\}.$$

Since $J\not\in \mathcal{C}_{0}$, we have that (\ref{defn:c0}) is 
false. It follows that $D_{J}\neq\emptyset$, so that we can define
$c_{J}:=a_{J}$ and $d_{J}:=\inf D_{J}<b$. 

Roughly speaking, if $(a,b)\setminus\Omega_{1}$ is a finite set, then
$[c_{J},d_{J}]$ is the closure of the union of $(a_{J},b_{J})$ and of
the first $2k$ connected components of $\Omega_{1}$ on its right.  If
$(a,b)\setminus\Omega_{1}$ is not finite, then we stop the union as
soon as we find an accumulation point of $(a,b)\setminus\Omega_{1}$.
For example, $d_{J}=b_{J}$ when $b_{J}$ is an accumulation point of
$(a,b)\setminus\Omega_{1}$.

We have to show that properties~(P1) through~(P3) are satisfied.  This
is trivial for~(P1) because $c_{J}=a_{J}>a$ (we recall that the
connected component of $\Omega_{1}$ with $a_{J}=a$, if it exists, does
not belong to $\mathcal{C}_{1}$), and $b_{J}\leq d_{J}<b$.

In order to prove~(P2), let us assume by contradiction that
(\ref{eqn:P2}) is false for some $x_{0}\in(a,b)$.  Then
$x_{0}\in(c_{J},d_{J})$ for at least $(2k+2)$ connected components
$J_{0}$, $J_{1}$, \ldots, $J_{2k+1}$ in $\mathcal{C}_{1}$.  Let us
assume, without loss of generality, that components are named in such
a way that $b_{J_{0}}<b_{J_{1}}<\ldots<b_{J_{2k+1}}$.

On the one hand we have
that $[b_{J_{0}},b_{J_{2k}}]\setminus\Omega_{1}\supseteq\{
b_{J_{0}},b_{J_{1}},\ldots,b_{J_{2k}}\}$, hence this set has at least
$(2k+1)$ elements.  This means that $b_{J_{2k}}\in D_{J_{0}}$, hence
$b_{J_{2k}}\geq d_{J_{0}}>x_{0}$.  On the other hand we have also that
$x_{0}>c_{J_{2k+1}}=a_{J_{2k+1}}\geq b_{J_{2k}}$, which gives a
contradiction.

It remains to prove (P3).  Let us take any $x\in D_{J}$.  Then
$[b_{J},x]$ contains at least $(2k+1)$ points which are not in
$\Omega_{1}$, hence at least $(2k+1)$ points where either $g$ vanishes
or $g'$ vanishes.  As a consequence, $[b_{J},x]$ contains either at
least $(k+1)$ points where $g$ vanishes, or at least $k$ points where
$g'$ vanishes.  In the first case Rolle's Theorem implies the
existence of at least $k$ points where $g'$ vanishes, so actually we
are always in the second case.

Applying Rolle's Theorem again, we obtain that $[b_{J},x]$ contains at
least $(k-1)$ points where $g''$ vanishes, at least $(k-2)$ points
where $g'''$ vanishes, and so on, up to find at least one point where
$g^{(k)}$ vanishes.  We have thus proved that all derivatives of $g$
up to order $k$ vanish at least once in $[b_{J},x]$.  Since they are
continuous functions, the same is true in $[b_{J},d_{J}]$, hence also
in $[c_{J},d_{J}]$.  This completes the proof of~(P3).

\subparagraph{\textmd{\emph{Absolute continuity in the union of
$\mathcal{C}_{1}$-components}}}

Let $J\in \mathcal{C}_{1}$, and let $[c_{J},d_{J}]$ be the interval
defined in the previous paragraph.  Due to (P1) and (P3) we can apply
Proposition~\ref{prop:cd} in the intervals
$(a_{J},b_{J})\subseteq[c_{J},d_{J}]$.  We obtain that 
$$\|f'\|^{p}_{p,w,J} \leq [C_{2}(k)]^{p}
\left[\hold\left(g^{(k)},(c_{J},d_{J})\right)\right]^{p/(k+\alpha)}
\cdot(d_{J}-c_{J}).$$

Estimating the H\"{o}lder constant in $(c_{J},d_{J})$ with the
H\"{o}lder constant in $(a,b)$, and summing over all $J\in
\mathcal{C}_{1}$, we deduce that $$\sum_{J\in \mathcal{C}_{1}}
\|f'\|^{p}_{p,w,J}\leq [C_{2}(k)]^{p}
\left[\hold\left(g^{(k)},(a,b)\right)\right]^{p/(k+\alpha)}
\cdot\sum_{J\in \mathcal{C}_{1}}(d_{J}-c_{J}).$$

Now we estimate the series in the right-hand side through a
``geometric double counting'' argument. Let $\chi_{(c_{J},d_{J})}(x)$ 
be 1 if $x\in(c_{J},d_{J})$, and 0 otherwise. Due to (\ref{eqn:P2}) we 
have that
$$\sum_{J\in \mathcal{C}_{1}}(d_{J}-c_{J})=
\sum_{J\in \mathcal{C}_{1}}\int_{a}^{b}\chi_{(c_{J},d_{J})}(x)\,dx=
\int_{a}^{b}\left(\sum_{J\in \mathcal{C}_{1}}
\chi_{(c_{J},d_{J})}(x)\right)\,dx\leq
(2k+1)(b-a),$$
hence
\begin{equation}
	\sum_{J\in \mathcal{C}_{1}} \|f'\|^{p}_{p,w,J}\leq
	(2k+1)[C_{2}(k)]^{p}
	\left[\hold\left(g^{(k)},(a,b)\right)\right]^{p/(k+\alpha)}
	\cdot(b-a).
	\label{est:c1}
\end{equation}

From (\ref{est:c0}) and (\ref{est:c1}) we easily obtain 
(\ref{th:main}).\qed

\subsubsection*{Proof of Theorem~\ref{thm:main-ndim}}

Let us assume, without loss of generality, that
$\Omega':=(a_{1},b_{1})\times\ldots\times(a_{n},b_{n})$
is a rectangle. From statement~(1) of Lemma~\ref{lemma:weak-lp} we 
have that
$$\|\nabla f\|_{p,w,\Omega'}\leq 
\left\|\sum_{m=1}^{n}|f_{x_{m}}|\right\|_{p,w,\Omega'}\leq
n\sum_{m=1}^{n}\|f_{x_{m}}\|_{p,w,\Omega'}.$$

Thus we can limit ourselves to
showing that all partial derivatives of $f$ lie in
$L^{p}_{w}(\Omega')$.  Without loss of generality, we prove this fact
in the case $m=1$.  Let $y:=(x_{2},\ldots,x_{n})$ denote the vector of
the remaining variables, and let
$\Omega'':=(a_{2},b_{2})\times\ldots\times(a_{n},b_{n})$.  Then for
every $M\geq 0$ we have that
\begin{equation}
	\mathrm{meas}\left\{x\in\Omega':|f_{x_{1}}(x)|>M\right\}=
	\int_{\Omega''}
	\mathrm{meas}\left\{x_{1}\in(a_{1},b_{1}):|f_{x_{1}}(x_{1},y)|>M\right\}
	\,dy.
	\label{fubini}
\end{equation}

For every fixed $y\in\Omega''$, the measure in the integral involves
only the function $x\to f(x,y)$, which is a function of one real
variable.  Therefore the measure can be estimated with the aid of
Theorem~\ref{thm:main} in terms of the function of one real variable
$x\to g(x,y)$, which for simplicity we denote by $g(\cdot,y)$.  We
obtain that 
\begin{eqnarray}
	\lefteqn{\hspace{-0.8em}
	\mathrm{meas}\left\{x_{1}\in(a_{1},b_{1}):
	|f_{x_{1}}(x_{1},y)|>M\right\} \ \leq\  
	\frac{1}{M^{p}}\,\|f_{x_{1}}(\cdot,y)\|_{p,w,(a_{1},b_{1})}^{p}} 
	\nonumber  \\
	\noalign{\vspace{1ex}}
	\hspace{0.9em} & \leq & \frac{[C(k)]^{p}}{M^{p}}\max\left\{
	\left[\hold\left(D^{k}_{x_{1}}g(\cdot,y),
	(a_{1},b_{1})\right)\right]^{p/(k+\alpha)}
	(b_{1}-a_{1}),
	\|g_{x_{1}}(\cdot,y)\|^{p/(k+\alpha)}_{L^{\infty}((a_{1},b_{1}))}\right\}
	\nonumber \\
	 \noalign{\vspace{1ex}}
	 & \leq & \frac{\left[C(k)\right]^{p}}{M^{p}}\max\left\{
	 \left[\hold\left(g^{(k)},\Omega\right)\right]^{p/(k+\alpha)}
	 \cdot(b_{1}-a_{1}), \|\nabla
	 g\|^{p/(k+\alpha)}_{L^{\infty}(\Omega)}\right\}.
	 \label{est:meas}
\end{eqnarray}

Therefore we can estimate the left-hand side of (\ref{fubini}) with
the right-hand side of (\ref{est:meas}) multiplied by the
$(n-1)$-dimensional measure of $\Omega''$.

Thus we obtain that
\begin{eqnarray*}
	\|f_{x_{1}}\|_{p,w,\Omega'} & = & \sup_{M\geq 0}M\cdot
	\left[\mathrm{meas}
	\left\{x\in\Omega':|f_{x_{1}}(x)|>M\right\}\right]^{1/p}
	\\
	 & \leq & C(k)\max\left\{
	 \left[\hold\left(g^{(k)},\Omega\right)\right]^{1/(k+\alpha)}
	 \cdot|\Omega'|^{1/p}, \|\nabla
	 g\|^{p/(k+\alpha)}_{L^{\infty}(\Omega)}
	 \cdot|\Omega''|^{1/p}\right\},
\end{eqnarray*}
where $|\Omega'|$ and $|\Omega''|$ are the $n$-dimensional measure of 
$\Omega'$, and the $(n-1)$-dimensional measure of $\Omega''$, 
respectively.

This is enough to complete the proof.\qed

\setcounter{equation}{0}
\section{Examples and open problem}\label{sec:ex}

\begin{ex}\label{ex:glaeser-v'}
	\begin{em}
		Let us consider the function $v(x):=\sin^{2}x+e^{-x^{2}}$.
		Then $v\in C^{\infty}(\re)$, and all its derivatives are
		bounded in $\re$, which means that $v\in C^{k,\alpha}(\re)$
		for every admissible value of $k$ and $\alpha$.  Moreover
		$v(x)>0$ for every $x\in\re$.  Let us consider the sequence
		$x_{n}:=2\pi n+n^{-1}$.  It is not difficult to see that
		$$v(x_{n})\sim\frac{1}{n^{2}}, \quad\quad v'(x)\sim\frac{2}{n}
		\quad\quad
		\mbox{as }n\to +\infty,$$
		hence the ratio (\ref{ratio}) is not bounded whenever $k\geq
		2$ and $\alpha>0$.
		
		This example shows that the conclusion of 
		Theorem~\ref{thm:glaeser-k} is false when $v'(x)$ is a sign 
		changing function (as in this example).
	\end{em}
\end{ex}

\begin{ex}\label{ex:no-glaeser}
	\begin{em}
		Let $\varphi\in C^{\infty}(\re)$ be a function such that 
		$\varphi(x)=1$ for every $x\leq 0$, $\varphi(x)=0$ for every 
		$x\geq 1$, and $\varphi'(x)<0$ for every $x\in(0,1)$. We note 
		that these conditions imply that all derivatives 
		of $\varphi$ vanish in $x=0$ and $x=1$. Let us consider the 
		following sequences
		$$\gamma_{n}:=e^{-(n+1)^{2}},
		\hspace{3em}
		a_{n}:=\gamma_{n}^{k+\alpha}\cdot|\log\gamma_{n}|,
		\hspace{3em}
		\alpha_{n}:=\sum_{i=n+1}^{\infty}a_{i}.$$
		
		Some simple calculus shows that $\alpha_{n}$ is well defined 
		(namely the series converges), and
		\begin{equation}
			\lim_{n\to+\infty}\frac{\alpha_{n}}{a_{n}}=0.
			\label{ex:limite}
		\end{equation}
		
		Let $w:\re\to\re$ be defined by
		$$w(x):=\left\{
		\begin{array}{ll}
			\alpha_{0}+a_{0} & \mbox{if }x\leq 0,  \\
			\alpha_{n}+a_{n}\varphi\left(\gamma_{n}^{-1}(x-n)\right)\quad & 
			\mbox{if }x\in[n,n+\gamma_{n}],  \\
			\alpha_{n} & \mbox{if }x\in[n+\gamma_{n},n+1].
		\end{array}
		\right.$$
		
		It is easy to see that $w\in C^{\infty}(\re)$, and for every 
		$x\in\re$ we have that $w(x)>0$ and $w'(x)\leq 0$. Moreover 
		for every $\beta\in(0,1]$ we have that
		$$w\in C^{k,\beta}(\re)\Longleftrightarrow
		\sup_{n\in\n}\frac{a_{n}}{\gamma_{n}^{k+\beta}} <+\infty,$$
		hence $w\in C^{k,\beta}(\re)$ for every $\beta<\alpha$, but not 
		for $\beta=\alpha$.
		
		Let us consider now the function $\psi(x):=\arctan
		e^{-(x-\gamma_{0}/2)^{3}}$, and let us set
		$v(x):=w(x)+\psi(x)$.  It is not difficult to see that
		$\psi\in C^{\infty}(\re)$, its derivatives of any order are
		bounded in $\re$, and $\psi(x)>0$, hence also $v(x)>0$, for
		every $x\in\re$.  Moreover we have that $\psi'(x)\leq 0$ for
		every $x\in\re$, with equality if and only if $x=\gamma_{0}/2$
		(where $w'(x)$ does not vanish), so that $v'(x)<0$ for every
		$x\in\re$.  Finally we have that $v\in C^{k,\beta}(\re)$ for every
		$\beta<\alpha$, but not for $\beta=\alpha$.
		
		Let us consider now the sequence $z_{n}:=n+\gamma_{n}/2$. It 
		is not difficult to check that
		$$\lim_{n\to +\infty}\frac{\psi(z_{n})}{a_{n}}=0,
		\hspace{3em}
		\lim_{n\to 
		+\infty}\frac{\psi'(z_{n})\cdot\gamma_{n}}{a_{n}}=0.$$
		
		Exploiting also (\ref{ex:limite}), it follows that (up to
		numeric constants) we have that
		$$v(z_{n})=a_{n}\varphi(1/2)+\alpha_{n}+\psi(z_{n})
		\sim a_{n}
		\quad\mbox{as }n\to +\infty,$$
		$$v'(z_{n})=\frac{a_{n}}{\gamma_{n}}\,\varphi'(1/2)+\psi'(z_{n})
		\sim \frac{a_{n}}{\gamma_{n}}
		\quad\mbox{as }n\to +\infty.$$
		
		This implies that the sequence
		$$\frac{|v'(z_{n})|^{k+\alpha}}{[v(z_{n})]^{k+\alpha-1}}\sim
		\frac{a_{n}^{k+\alpha}}{\gamma_{n}^{k+\alpha}}\cdot
		\frac{1}{a_{n}^{k+\alpha-1}}=
		\frac{a_{n}}{\gamma_{n}^{k+\alpha}}=
		|\log\gamma_{n}|$$
		is unbounded, hence the ratio (\ref{ratio}) is unbounded. 
	\end{em}
\end{ex}

\begin{ex}\label{ex:p}
	\begin{em}
		Let $f:(-1,1)\to\re$ and $g:(-1,1)\to\re$ be defined by
		$$f(x):=|x|^{1/(k+\alpha)},
		\quad\quad
		g(x):=x.$$
		
		All assumptions of Theorem~\ref{thm:main} are satisfied (and
		$g$ is actually of class $C^{\infty}$), and $f'\not\in
		L^{p}((-1,1))$ when $p$ is defined by (\ref{defn:p*}).  This
		shows that the summability of $f'$ given by
		Theorem~\ref{thm:main} is optimal.
	\end{em}
\end{ex}

\begin{ex}\label{ex:holder}
	\begin{em}
		Let $\varphi\in C^{\infty}(\re)$ be a function such that
		$\varphi(x)>0$ for every $x\in(0,1)$, and $\varphi(x)=0$
		elsewhere, so that all its derivatives vanish in $x=0$ and
		$x=1$.  Let us assume also that $\varphi'(x)>0$ for every
		$x\in(0,1/2)$, and $\varphi'(x)<0$ for every $x\in(1/2,1)$, so
		that it is quite easy to estimate the total variation of
		$\varphi$ and all its roots.  Let us consider the following
		sequences
		$$\gamma_{n}:=\frac{1}{(n+3)\log(n+3)\log^{2}(\log(n+3))},
		\quad\quad \lambda_{n}:=\sum_{i=n}^{\infty}\gamma_{i},
		\quad\quad
		a_{n}:=\gamma_{n}^{k+\alpha}\cdot|\log\gamma_{n}|.$$
		
		Some simple calculus shows that $\lambda_{n}$ is well defined 
		(namely the series converges). Let $g:(0,\lambda_{0})\to\re$ 
		be defined by
		$$g(x):=a_{n}\varphi\left(
		\frac{x-\lambda_{n+1}}{\lambda_{n}-\lambda_{n+1}}\right)
		\quad\quad
		\forall x\in[\lambda_{n+1},\lambda_{n}),$$
		and let $f(x):=[g(x)]^{1/(k+\alpha)}$ be its $(k+\alpha)$-th 
		root. 
		
		It is easy to see that $g\in 
		C^{\infty}((0,\lambda_{0}))$. Moreover 
		for every $\beta\in(0,1]$ we have that
		$$g\in C^{k,\beta}((0,\lambda_{0}))\Longleftrightarrow
		\sup_{n\in\n}\frac{a_{n}}{\gamma_{n}^{k+\beta}} <+\infty,$$
		hence $g\in C^{k,\beta}((0,\lambda_{0}))$ for every
		$\beta<\alpha$, but not for $\beta=\alpha$.  Moreover the
		total variation of $f$ in $(0,\lambda_{0})$ is given by
		$$TV(f,(0,\lambda_{0}))=2\sum_{n=0}^{\infty}a_{n}^{1/(k+\alpha)}
		\cdot \left[\varphi(1/2)\right]^{1/(k+\alpha)},$$
		hence
		$$f'\in L^{1}((0,\lambda_{0}))\Longrightarrow f\in
		BV((0,\lambda_{0}))\Longrightarrow
		\sum_{n=0}^{\infty}a_{n}^{1/(k+\alpha)}<+\infty$$
		(actually also the reverse implications hold true, but we 
		don't need them). Now some simple calculus shows that
		$$\sum_{n=0}^{\infty}a_{n}^{1/(k+\alpha)}=
		\sum_{n=0}^{\infty}\gamma_{n}|\log\gamma_{n}|^{1/(k+\alpha)}=
		+\infty.$$
			
		This means that $f(x)$ is not a bounded variation function in
		$(0,\lambda_{0})$, and hence $f(x)$ is not even absolutely
		continuous in $(0,\lambda_{0})$.
	\end{em}
\end{ex}

We conclude by recalling the main open problem in this field, namely
the absolute continuity of roots of polynomials instead of roots of
functions.  We point out that the problem is open even in the case of
a real root of a non-hyperbolic polynomial with real coefficients of
class $C^{\infty}$.

\paragraph{Open problem}

Let $k$ be a positive integer, let $(a,b)\subseteq\re$ be an 
interval, and let
$$P(x):=x^{k+1}+a_{k}(t)x^{k}+\ldots+a_{1}(t)x+a_{0}(t)$$
be a monic polynomial with complex valued coefficients 
$a_{i}\in C^{k,1}((a,b))$. Let $z:(a,b)\to\C$ be a continuous 
function such that $P(z(t))=0$ for every $t\in(a,b)$.

Then $z'\in L^{p}_{w}((a,b))$ for $p=1+1/k$.

\subsubsection*{\centering Acknowledgments}

We discovered this problem thanks to a conference held by Prof.\ F.\
Colombini during the meeting ``Asymptotic Properties of Solutions to
Hyperbolic Equations'' (London, March 2011).  For this reason, we are
indebted to him and to the organizers of that workshop.  

We would like also to thank Prof.\ S.\ Spagnolo for pointing out
reference~\cite{tarama}, and Prof.\ F.\ Colombini for pointing out
reference~\cite{colombini-lerner}.

Finally, we thank the anonymous referee for carefully checking the
manuscript, and for many insightful and powerful suggestions which 
greatly improved the clarity of the presentation.

\label{NumeroPagine}

\end{document}